\documentclass[12pt]{article}
\usepackage{amsmath,amssymb,amsthm}

\usepackage{fullpage}

\title{On a theorem concerning the elliptic functions}

\author{EDVARD PHRAGMEN \\[12pt]
        Stolkolm\\
        Acta Mathematica, vol 7\\
        pp. 33-42\footnote{(translator's note) The following is a literal word-for-word translation from the French of the work "Sur un Theoreme concernant les fonctions elliptiques" by Edvard Phragmen.}}

\date{\today}

\theoremstyle{definition}

\begin{document}
\maketitle
In his lectures at the university of Berlin, \textsc{Weierstrass} has given a theory of elliptic functions notable for its beauty as well as for the simplicity of its methods, in which he takes the following theorem as his point of departure:

Every analytic function $\phi(u)$ possessing an addition theorem, in other words, such that there exists an algebraic relation between the values of the function corresponding to the values $u$, $v$, $u+v$, of the argument is:\begin{enumerate}
  \item either an algebraic function of $u$;
  
  \item or, if $\omega$ denotes a conveniently chosen constant, an algebraic function of the function $e^{\frac{\pi i u }{\omega}}$;
  \item or, if $\omega, \omega'$ denote two conveniently chosen constants, an algebraic function of the function\footnote{For the theory of this function, see the work of \textsc{Schwarz}: \emph{Formeln und Lehrs$\ddot{a}$tze zum Gebrauche der elliptischn Functionen } 1881-1884.}
  \begin{equation*}
\label{ }
\wp(u|\omega,\omega')=\frac{1}{u^2}+\sum' \left[\frac{1}{(u-2\mu \omega-2\mu'\omega')^2}-\frac{1}{(2\mu \omega+2\mu'\omega')^2} \right] \end{equation*}($\mu,\mu'=0,1,2\cdots$ excluding the combination $\mu=\mu'=0$.)
\end{enumerate}

This theory has never been published, in its entirety, by the author; the only work, to my knowledge, which gives some idea of it is that of \textsc{Schwarz}; however, since the demonstration which \textsc{Weierstrass} has given of the stated theorem has never been published, I have been reduced to my own proper resources to find one; only later was \textsc{Weierstras}' demonstration communicated to me by my professors, \textsc{Mittag-Leffler} and madmoselle \textsc{Kowalevski}.  An essential characteristic of this demonstration is that it may be generalized to the general case of the abelian functions.  However, if one renounces this genralization and limits oneself to the demonstration of the theorem just stated, then the demonstration which I have found may be a bit simpler and easier than that of \textsc{Weierstrass}, and I have the honor of presenting it to the readers of the \textsc{Acta} in the hope that it interests them.  
\\
\\
\\
\\
\\

We now assume that the analytic function $\phi(u)$ has an algebraic addition theorem in which an algebraic relation exists between $\phi(u)$, $\phi(v)$, and $\phi(u+v)$.  I limit myself to supposing that there are three elements \footnote{$\mathrm{P(u|a)}=$ an ordinary power series in $(u-a)$}\begin{equation*}
\label{ }
\mathrm{P_1(u|a)},\\ \mathrm{P_2(v|b)}, \\ \\ \mathrm{P_3(u+v|a+b)}
\end{equation*}of the function $\phi(u)$ connected by an algebraic equation\begin{equation}\label{AT}
G\left[\mathrm{P_1(u|a)},\\ \mathrm{P_2(v|b)}, \\ \\ \mathrm{P_3(u+v|a+b)}\right]=0
\end{equation}It is well known that the same equation continues to hold for every system of three elements which may be obtained by continuing the system of elements $\mathrm{P_1, P_2, P_3}$ along any path whatsoever.  However, we may not conclude immediately that one has\begin{equation*}
\label{ }
G\left[\phi(u),\phi(v),\phi(u+v)\right]=0
\end{equation*}for any three values whatsoever of the function corresponding to the values $u$, $v$, $u+v$, of the argument.  In fact, if we change the element $\mathrm{P_1(u|a)}$ into $\mathrm{P_1'(u|a)}$ by continuing it along a closed path while leaving the element $\mathrm{P_2(v|b)}$ unaltered, generally the element $ \mathrm{P(u+v|a+b)}$ will change into a new element  $\mathrm{P'(u+v|a+b)}$ such that one does not know whether or not the new element $\mathrm{P_1'(u|a)}$ satisfies the relation \begin{equation*}
\label{ }
G\left[\mathrm{P'_1(u|a)},\\ \mathrm{P_2(v|b)}, \\ \\ \mathrm{P(u+v|a+b)}\right]=0.
\end{equation*}
\\

In the following I will prove, in the first place, that if one admits the relation \eqref{AT}, the function $\phi(u)$ has the character of an algebraic function in the interior of any finite domain.

In saying that a function defined by a given element has the character of an algebraic function in the interior of a given domain, I mean that all of its values which may be obtained by continuing the element along paths which totally belong to the given domain satisfy an equation of the form\begin{equation*}
\label{ }
z^n+f_1z^{n-1}+f_2z^{n-2}+\cdots+f_n=0
\end{equation*}where $f_1$, $f_2$, $\cdots$, $f_n$ have the character of a rational function in the interior of the same domain.

Then I will show that one has in effect\begin{equation*}
\label{ }
G\left[\phi(u),\phi(v),\phi(u+v)\right]=0
\end{equation*}where $\phi(u)$ denotes any value whatsoever of the function $\phi$ corresponding to the argument $u$, and similarly for  $\phi(v)$ and  $\phi(u+v)$.

Consequently, I know that  $\phi(u)$ can only have a finite number of values at each point, and based on this fact, I will prove that  $\phi(u)$ is an algebraic function of a periodic function, and there is little difficulty in completing the proof
\\
\\
We suppose, therefore, that the relation\begin{equation*}\label{AAT}
G\left[\mathrm{P_1(u|a)},\\ \mathrm{P_2(v|b)}, \\ \\ \mathrm{P(u+v|a+b)}\right]=0
\end{equation*}holds between three elements $\mathrm{P_1}$, $\mathrm{P_2}$, and $\mathrm{P}$ of the function $\phi$.

If we consideer one of these elements, for example $\mathrm{P_1(u|a)}$, we see that one may take a positive quantity $r_1$ such that the function defined by this element has the character of an algebraic function within the domain\begin{equation*}
\label{ }
|u-a|<r_1.
\end{equation*}This all the more true if we choose $r_1$  smaller than the radius of convergence of the element $\mathrm{P_1(u|a)}$.  Analogous results hold for the elements $\mathrm{P_2(v|b)}$ and $\mathrm{P(u+v|a+b)}$.

Moreover, it is seen that the upper limits of the quantities $r_1$, $r_2$, $r$, are all finite or infinite at the same time.  For, if one of the upper limits be infinite, that is that $\phi(u)$ has the character of an algebraic function for every finite domain, then it is necessary that the others also be infinite.   But, it is easy to see that all three cannot be finite.  For, if this be the case, let $\rho_1$, $\rho_2$, $\rho$, be their values: then it is necessary that one of the following two inequalities hold:$$\rho_1<\rho +\rho_2 \ \ \text{or}\ \ \rho_2<\rho+\rho_1.$$  Suppose $\rho_1<\rho +\rho_2 $.  Then, if we put $$v-b=-\frac{\rho_2}{\rho+\rho_2}(u-a),$$and therefore $$u+v-a-b=\frac{\rho}{\rho+\rho_2}(u-a),$$and if we suppose$$|u-a|<R<\rho+\rho_2,$$one has $$|v-b|<\frac{R}{\rho+\rho_2}\rho_2<\rho_2,$$ $$|u+v-a-b|<\frac{R}{\rho+\rho_2}\rho<\rho,$$and consequently, the functions of $u$ defined after this substitution by the elements $\mathrm{P_2}$ and $\mathrm{P}$ have the character of an algebraic function in all of the domain $$|u-a|<R,$$where $R$ is a positive quantity smaller than $\rho+\rho_2.$  If one denotes these functions by $y,z$ and the function defined by the element $\mathrm{P_1}$ by $x$, we have the relation $$G(x,y,z)=0.$$  If we eliminate $y$ and $z$ between this equation and the equations $$y^{\mu}+\phi_1y^{\mu-1}+\cdots+\phi_{\mu}=0,$$ $$z^{\nu}+\psi_1y^{\nu-1}+\cdots+\psi_{\nu}=0,$$which define $y$ and $z$ as functions of algebraic character in the interior of the domain $|u-a|<R$, we obtain an equation $$x^n+f_1x^{n-1}+\cdots+f_n=0,$$where $f_1$, $\cdots$, $f_n$ have the character of a rational function in the domain $|u-a|<R$.

Therefore, $x$ has the character of an algebraic function in the interior of this same domain $|u-a|<R$, where $R$ is any positive quantity whatsoever smaller than $\rho+\rho_2,$ which contradicts the supposition $\rho_1<\rho+\rho_2.$

We see that the same reasoning applies to the case $\rho_2<\rho+\rho_1.$

Therefore $\phi(u)$ has the character of an algebraic function in every finite domain.

We now prove that we always have the relation \begin{equation*}
\label{ }
G\left[\phi(u),\phi(v),\phi(u+v)\right]=0
\end{equation*} between three values of the function $\phi$ corresponding to the values $u,v,u+v$ of the argument.

We will prove this by showing that each value $\phi_{a}$, which the function $\phi(u)$ can take at the point  $a$, which $\phi(u)$ approaches indefinitely when  $u$ approaches $a$ indefinitely along a certain path, satisfies the equation\begin{equation*}\label{AAT}
G\left[\phi_a,\\ \mathrm{P_2(v|b)}, \\ \\ \mathrm{P(u+v|a+b)}\right]=0.
\end{equation*}

In fact, one arrives at the vaue $\phi_a$ by continuing the element $\mathrm{P_1(u|a)}$ along a certain closed path.  One may always chose a continuous, finite and simply connected domain which contains this path completely.  One may then detrmine a second continuous and finite domain such that the lower limit of the distances of an interior point of the first domain and a point in the exterior of the second are larger than a quantity $d>|b|$.
We imagine, for a moment, a third domain large enough to contain in its interior not only the last of the two domains we have just fixed, and its boundary, but also the path along which it is necessary to continue the element $\mathrm{P_1(u|a)}$ in order to arrive at the element $\mathrm{P(u+v|a+b)}$.  In the interior of this domain, the function defined by the element $\mathrm{P_1(u|a)}$--or, which gives the same function, by $\mathrm{P(u+v|a+b)}$--has the character of an algebraic function, and therefore it has a finite number of singular points in the interior of the larger of the first two domains.  We may therefore choose two points, $a'$, $b'$ in the neighborhood of $a$ and $b$ such that they satisfy the following conditions: in the first place it is necessary that we have $|b'|<d$,   in order that to each point $\alpha$ in the interior of the smaller of the two domains corresponds a point $\alpha+b'$ in the interior of the larger.  Moreover, the point $a'$ will be chosen in the interior of the small domain in such a way that the points $a'$, $a'+b'$ do not coincide with any of th singular points, of which there are a finite number in the large domain which constitute the singular points of the function considered up to now.  Finally, if $\alpha$ is a singular point of this function situated in the interior of the small domain, $\alpha+b'$ will be a regular point of the same function.

Without changing the final value $\phi_a$, we may now replace the closed path $a\cdots a$ by a set of regular paths as follows: $1^{\text{o}}$ a path $aa'$, $2^{\text{o}}$ a number of \emph{loops}\footnote{By a \emph{loop} joining the point $a'$ to a singular point $a$, we mean a regular path composed of a path leaving $a'$ and arriving at the neighborhood of $a$, a small circle around $a$, and the first path described in the opposite direction .  (See, for ex, \textsc{Briot} and \textsc{Bouquet, \emph{Th\'eorie des fonctions elliptiques).}}} joining the point $a'$ to the singular points situated in the smaller domain, and $3^{\text{o}}$ the path $a'a$.  We may choose these loops in such a way that the corresponding loops which the point $u+b'$ describes are composed of regular paths and of circles which do not contain any singular point.  Therefore, if one continues the system of three elements $\mathrm{P_1(u|a)},\mathrm{P_2(v|b)},  \mathrm{P(u+v|a+b)}$ by varying $1^{\text{o}}$ $u$ from $a$ to $a'$ and $v$ from $b$ to $b'$, $2^{\text{o}}$ $u$ from $a'$ to $a'$ along the loops, and $3^{\text{o}}$ $u$ from $a'$ to $a$ and $v$ from $b'$ to $b$, by starting with the element $\mathrm{P_1(u|a)}$ one approaches indefinitely the value $\phi_a$, while the elements $\mathrm{P_2(v|b)}$ and $\mathrm{P(u+v|a+b)}$ return to themselves.

But, it has therefore been demonstrated that every regular or singular element $\phi(u|a)$ of the function $\phi(u)$ in the neighborhood of the point $a$ necessarily satisfies the equation\begin{equation*}\label{AAT}
G\left[\phi(u|a),\mathrm{P_2(v|b)},  \mathrm{P(u+v|a+b)}\right]=0.
\end{equation*} Since the analogous result holds for the elements $\mathrm{P_2}$ and $\mathrm{P}$, and since moreover the points $a,b$ were taken totally arbitrarily, our assertion is completely justified.

From this fact that the relation\begin{equation*}
\label{ }
G\left[\phi(u),\phi(v),\phi(u+v)\right]=0
\end{equation*}always holds it follows that the function $\phi(u)$ has only a finite number of values at any point $$\phi_1(u),\phi_2(u), \cdots, \phi_n(u).$$Therefore, if we form the expressions

\begin{align*}
\frac{1}{\phi_1u-a}&+\cdots+\frac{1}{\phi_nu-a}  \\
\frac{1}{(\phi_1u-a)^2}&+\cdots+\frac{1}{(\phi_nu-a)^2} \\
...............&..............................\\
\frac{1}{(\phi_1u-a)^n}&+\cdots+\frac{1}{(\phi_nu-a)^n}
\end{align*}we see immediately that they represent functions which have the character of a rational function in the interior of all finite domains.  If they are all rational functions, $\phi(u)$ is an algebraic function.  Otherwise, suppose that the function$$\frac{1}{(\phi_1u-a)^m}+\cdots+\frac{1}{(\phi_nu-a)^m}$$has the point $u=\infty$ as an essential singular point.  Then we know that for  values of the absolute value of the argument larger than any given quantity, this function may take values larger in absolute value than any given quantity, and that therefore for values of $u$ larger in absollute value than  any given quantity, $\phi(u)$ may take values such that the absolute value of $[\phi(u)-a]$ may be smaller than any given quantity.  But, this being the case, we demonstrate with the help of reasoning employed by \textsc{Weierstrass} in his course and which I reproduce here succinctly, that $\phi(u)$ is necessarily periodic.

Indeed, we see in the first place that in the neighborhood of a given point, there is always a value which $\phi(u)$ takes on at least any given number of times.  For, suppose that in the neighborhood of a point $a$ there is a value $b$ which $\phi(u)$ takes on $m$ times at the points $u_1,\cdots,u_m$.  But, $\phi(u)$ may take every value in a certain neighborhood of $b$ for the values of the argument in the neighborhoods of the points $u_1,\cdots,u_m$.  However $\phi(u)$ may take a value as little different from $b$ as we like in a point in the neighborhood of $u=\infty$, so we see that we may take in each neighborhood of $b$ a value which $\phi(u)$ takes at least $(m+1)$ times.  In the neighborhood of each value $a$, one may therefore take another  value $b$ which $\phi(u)$ takes at least a given number of times.

Therefore, if the equation \begin{equation*}
\label{ }
G\left[\phi(u),\phi(v),\phi(u+v)\right]=0\end{equation*}is of degree $m$ with respect to $\phi(u+v)$ and if we choose $(m+1)$ points $v_1,\cdots,v_{m+1}$, in which  $\phi(v)$ takes the same value $a$, the equation in $z$ $$G[\phi(u), a, z]=0$$has the roots \newpage$$\phi_1(u+v_1), \cdots, \phi_1(u+v_{m+1})$$ $$.............................................$$  $$\phi_n(u+v_1), \cdots, \phi_n(u+v_{m+1}).$$

Therefore,  since the number of different roots may not exceed $m$, we have for each given value of $u$ an equality of the form  $$\phi_{\alpha}(u+v_{\gamma})=\phi_{\beta}(u+v_{\delta})\ \ \ \ \ \ \ \ \ \ \ \ \ \ \ (\gamma><\delta)$$from which, the number of these combinations being finite, at least one of these equalities must be an identity.  But, from the identity $$\phi_{\alpha}(u+2\omega)=\phi_{\beta}(u)$$ that is to say, the identity of one of the elements of $\phi(u+2\omega)$ at any point whatever with that of one of the elements of $\phi(u)$ at the same point, we conclude that the functions of $u$, namely $\phi(u+2\omega)$ and $\phi(u)$ are identical, that is to say, that $\phi(u)$ has the period $2\omega$.

Suppose now $2\omega$ is a period of $\phi(u)$ such that there are no periods of the form $2\mu \omega$, $\mu=\text{a real quantity smaller than unity}$.  Such a period exists.  Moreover, let $z=e^{\frac{\pi i u}{\omega}}$.  The elementary symmetric functions of  $\phi_1(u),\phi_2(u), \cdots, \phi_n(u)$ are uniform analytic functions of $z$ which can have at most the two singular points $z=0$ and $z=\infty$, and we see, as above, that if the function $\phi(u)$ is not algebraic in $z$, it takes on the same value for $(m+1)$ values of $z$, and therefore for $(m+1)$ non-equivalent values of $u$ with respect to the period $2\omega$.  From this we conclude that $\phi(u)$ has a period whose ratio with $2\omega$ is not real.

We now choose any two periods $2\omega, 2\omega'$ whose ratio is not a real quantity and we form the function $\wp(u|\omega,\omega')$ belonging to these periods: by putting $z=\wp(u|\omega,\omega')$ we easily see that $\phi(u)$ is an algebraic function of $z$.  For, from the definition of $\wp(u|\omega,\omega')$ it follows immediately that this function is a doubly periodic function, having elementary periods $2\omega,2\omega'$ which has no essential singularity at any finite distance.  Consequently, since $\wp(u)$ has a double infinity in the parallelogram of periods, we conclude that in each elementary parallelogram there are two distinct points or two points which coincide, where this function acquires any given value.  And, whenever $z$ is situated in the neighborhood of any point $z_0$, finite or infinite, each of the two corresponding values of $u$ which are non equivalent with respect to the periods have the character of an algebraic function of $z$.  Therefore, $\wp(u)$ considered as a function of $z$ also has the character of an algebraic function in the neighborhood of every point $z_0$ (finite or infinite), and therefore this function can only have a finite number of values for every value of $z$.  

But, from this we may conclude that $\phi(u)$ is an algebraic function of $z$.
\\
\\
\\
\\
\hrule\hrule

\end{document}